\documentclass[reqno]{amsart}
\usepackage{amscd}
\newtheorem{teo}{Theorem}[section]

\theoremstyle{definition}
\newtheorem{define}[teo]{Definition}

\newtheorem{lema}[teo]{Lemma}
\newtheorem{cor}[teo]{Corollary}
\newtheorem{prop}[teo]{Proposition}
\theoremstyle{remark}

\numberwithin{equation}{section}

\begin{document}
\newcommand{\slq}{\ensuremath{(sl_{n+1}^{+})_{q}\;}}
\newcommand{\an}{\ensuremath{\mathbb{A}_{n}}}
\newcommand{\Uq}{\ensuremath{U_{q}^{+}(sl_{n+1})\;}}
\newcommand{\B}{\ensuremath{\mathcal{B}\;}}
\newcommand{\s}{\ensuremath{\sigma}}
\newcommand{\xt}[1]{(x_{j_{#1}}^{n_{#1}-i_{#1}}\otimes x_{j_{#1}}^{i_{#1}})}
\newcommand{\xj}[1]{x_{j_{#1}}^{n_{#1}-i_{#1}}}
\newcommand{\xij}[1]{x_{j_{#1}}^{i_{#1}}}
\newcommand{\qbin}[1]{\begin{pmatrix}n_{#1}\\i_{#1}\end{pmatrix}_{q^{2}}}
\newcommand{\xnj}[1]{x_{j_{#1}}^{n_{#1}}}
\newcommand{\xmk}[1]{x_{k_{#1}}^{m_{#1}}}

\title[The positive part QUEA of type $\mathbb{A}_{n}$ as a braided quantum group]
{The positive part of the quantized universal enveloping algebra
of type \an as a braided quantum group}
\author[C. Bautista]{C\'esar Bautista}
\address{Facultad de Ciencias de la Computaci\'on.
Benem\'erita Universidad Aut\'onoma de Puebla. edif. 135. 14 sur y av.
San Claudio. Ciudad Universitaria. Puebla, Pue. c.p. 72570. M\'exico.}
\email{bautista@fcfm.buap.mx, bautista@solarium.cs.buap.mx}
\thanks{Partially supported by a grant of Fundaci\'on TELMEX
and by the Investigation Proyect INI106897  of DGAPA/UNAM}
\begin{abstract}
A generalized Hopf algebra structure for the positive (negative) part of
the Drinfeld-Jimbo quantum group of type \an is established without
make any use of the usual deformation of the abelian part of $sl_{n+1}$.
\end{abstract}
\maketitle

\begin{section}{Introduction}
The aim of this paper is to explain the bialgebra structure of the
positive part of the quantized universal enveloping algebra
(Drinfeld-Jimbo quantum group) of type \an using the Lie
algebra theory concepts.

Recently has been introduced a generalization of Lie algebras, the basic
$T$-Lie algebras \cite{yo}. Using the $T$-Lie algebra concept some new
(we think) quantum groups of type \an can be constructed.
Such quantum groups arise as universal enveloping algebras of certain
deformations as generalized Lie algebras of the Lie algebras form by
upper triangular matrices $sl_{n+1}^{+}$.

Let us explain, embedded in the positive (negative) parts
$U_{q}(sl_{n+1})$ of the Drinfeld-Jimbo quantum groups of type \an there
are some generalized Lie algebras \slq called $T$-Lie algebras,
\cite{yo}. Such $T$-Lie algebras satisfy not only a generalized
antisymmetry  and a generalized Jacobi identity, but an additional
condition called multiplicativity. Through these $T$-Lie algebras the
Poincar\'e-Birkhoff-Witt theorem for $U_{q}^{\pm}(sl_{n+1})$ can be
explained.

The   Poincar\'e-Birkhoff-Witt theorem  is a general property for the
universal enveloping algebras of adequate $T$-Lie algebras. In order to
keep the proof of such theorem closest to the classical one
\cite{Bourbaki} a Gurevich's condition of multiplicativity is needed,
\cite{yo}, \cite{Gurevich}.

On the other hand, the next natural step in the study of \slq as a
$T$-Lie algebra is to give to its universal enveloping algebra a
structure of Hopf algebra.

Now a structure as a generalized Hopf algebra (braided quantum group) of
the universal enveloping algebra \Uq  of \slq is
presented. As a matter of fact, \Uq has the usual algebra structure:
generators $x_{1},\ldots, x_{n}$ and relations
\begin{gather}
x_{i}x_{j}-x_{j}x_{i}=0,\text{ if }|i-j|>1\\
x_{i}x_{j}^{2}-(q+q^{-1})x_{i}x_{j}x_{i}+x^{2}_{j}x_{i}=0,\text{ if
}|i-j|=1
\end{gather}
but now the tensorial product $\Uq\otimes\Uq$  has a non-standard
algebra structure: the multiplication  is given by
\[
(a\otimes b)(c\otimes d)=a\sigma(b\otimes c)d
\]
where $\sigma:\Uq\otimes\Uq\rightarrow \Uq\otimes \Uq$  is not the usual
switch. This is so, because \Uq has a coproduct given by
\[
\phi(x_{i})=x_{i}\otimes 1+1\otimes x_{i},\;i=1,\ldots,n.
\]
We obtain a non-trivial deformation of the classical
universal enveloping algebra $U(sl_{n+1}^{+})$ without make any use of
the usuals $K_{1}^{\pm},\ldots K_{n}^{\pm}$ (which they form a
deformation of the abelian part of $sl_{n+1}$). Such developments leads
to the following question: are there any non-trivial deformation of
$U_{q}(sl_{n+1})$  constructed without the usuals
$K_{1}^{\pm},\ldots K_{n}^{\pm}$ ? We do not have an answer yet.
\end{section}
\begin{section}{Braids and Coproduct}
Let $k$ be a unitary commutative ring and $q\in k^{*}$. Denote with $e_{ij},$
 $1\leq i,j\leq n$ canonical basis of $gl_{n}$
the matrices $n\times n$, with \B the canonical basis of $sl_{n}^{+}$,
and with $[,]$ the usual bracket in $gl_{n}$. Define $m=n(n+1)/2$ and
$c_{ij,ab}\in \mathbb{Z}$, $1\leq i<j\leq m, 1\leq a<b\leq m$  such that
$[[e_{ij},e_{ji}],e_{ab}]=c_{ij,ab}e_{ab}$. Put
\[
e_{ij}<e_{ab}\text{ if }i+j<a+b\text{ or }(i+j=a+b \text{ and }j<b).
\]
Let \slq be $sl_{n}^{+}$
with structure of $T$-Lie algebra. This means: \slq is $sl_{n}^{+}$ in
its structure of $k$-module, besides a bracket $[,]_{q}$ such that
\[
[e_{ij},e_{ab}]_{q}=[e_{ij},e_{ab}],\text{ if }e_{ij}\leq e_{ab},
\]
a $k$-linear morphism
$S:\slq\otimes_{k}\slq\rightarrow\slq\otimes_{k}\slq$ called presymmetry
such that
\[
[,]_{q}S=-[,]_{q}
\]
defined by
\begin{gather*}
S(e_{ij}\otimes e_{ab})=q^{c_{ij,ab}}e_{ab}\otimes e_{ij},\text{ if
}e_{ij}<e_{ab}\\
S(e_{ij}\otimes e_{ij})=e_{ij}\otimes e_{ij},\;S^{2}=1
\end{gather*}
and a $k$-linear morphism
$\langle,\rangle:\slq\otimes_{k}\slq\rightarrow \slq\otimes_{k}\slq$
called pseudobracket such that
\[
[,]_{q}\langle,\rangle=0,\langle,\rangle S=-\langle,\rangle,
\]
in the case $e_{ij}<e_{ab}$ defined by
\[
\langle e_{ij},e_{ab}\rangle=\begin{cases}
                              (q-q^{-1})e_{aj}\otimes e_{ib},&\text{ if
                              }a<j<b\text{ and }i<a<j,\\
                              0,&\text{otherwise.}
                             \end{cases}
\]
The algebra \slq satisfies generalized Lie algebra axioms, see
\cite{yo}.

On the other hand, since
$[e_{ij},e_{ab}]=\delta_{jk}e_{il}-\delta_{li}e_{kj}$ where $\delta$ is
the Kronecker delta, we get
\begin{equation}\label{ces}
c_{ij,ab}=-\delta_{b,i}+\delta_{b,j}+\delta_{i,a}-\delta_{j,a},\;\forall
e_{ij},e_{ab}\in \B
\end{equation}
We shall define a new symmetry \s. Let
$\s:\slq\otimes_{k}\slq\rightarrow \slq\otimes_{k}\slq$ such that
\[
\s (e_{ij}\otimes e_{ab})=q^{c_{ij,ab}}e_{ab}\otimes e_{ij},\;\forall
e_{ij},e_{ab}\in\B.
\]
\begin{prop}\label{pse}
Let $[,]_{q}$ be the bracket of \slq. The linear morphism \s\; satisfies
the multiplicativity conditions:
\begin{gather}
\s(1\otimes [,]_{q})=([,]_{q}\otimes 1)\s_{23}\s_{12},\\
\s([,]_{q}\otimes 1)=(1\otimes [,]_{q})\s_{12}\s_{23}
\end{gather}
\end{prop}
\begin{proof}
The Jacobi identity in $sl_{n}^{+}$ ensures that
\[
[[e_{ij},e_{ji}],[e_{ab},e_{uv}]]=(c_{ij,ab}+c_{ij,uv})[e_{ab},e_{uv}]
\]
this implies
\begin{multline}
\s(e_{ij\otimes
[e_{ab},e_{uv}]_{q}})=q^{c_{ij,ab}+c_{ij,uv}}[e_{ab},e_{uv}]_{q}\otimes
e_{ij}\\=([,]_{q}\otimes 1)\s_{23}\s_{12}(e_{ij}\otimes e_{ab}\otimes
e_{uv}).
\end{multline}
In a similar way we obtain
\[
\s([e_{ij},e_{ab}]_{q}\otimes e_{uv})=(1\otimes
[,]_{q})\s_{23}\s_{12}(e_{ij}\otimes e_{ab}\otimes e_{uv})
\]
\end{proof}

Now we may extend \s\, to \Uq.\; To do so, let us consider $\mathcal{T}$
as the $k$-tensorial algebra of \slq. Define
$\tilde{\s}:\mathcal{T}\otimes\mathcal{T}\rightarrow
\mathcal{T}\otimes_{k}\mathcal{T}$   by
\[
\tilde{\s}(z_{1}\ldots z_{n}\otimes y_{1}\ldots
y_{m})=\prod_{i,j}p_{z_{i},y_{j}}y_{1}\ldots y_{m}\otimes z_{1}\ldots
z_{m},\;\forall z_{i},y_{j}\in\B,
\]
where $\s(z_{i}\otimes y_{j})=p_{z_{i},y_{j}}y_{j}\otimes
z_{i},\;p_{z_{i},y_{j}}\in k.$

In particular, the objects made by finite tensor products of one
dimensional spaces $\langle x\rangle$, $x\in\B$ form a {\it braided
tensor category} with braiding $\tilde{\s}$, \cite[Chapter 11]{Street}
and the proposition \ref{pse} says that the bracket $[,]_{q}$ is
compatible with the
braiding $\tilde{\s}$.
\begin{lema}
Let $\langle,\rangle$ be the pseudobracket of \slq and $m$ the product
of the tensor algebra $\mathcal{T}$. The linear morphism \s\, satisfies
the multiplicativity conditions
\begin{gather}
\tilde{\s}(1\otimes m\langle,\rangle)=(m\langle,\rangle\otimes
1)\tilde{\s}_{23}\tilde{\s}_{12},\\
\tilde{\s}(m\langle,\rangle\otimes 1)=(1\otimes
m\langle,\rangle)\tilde{\s}_{12}\tilde{\s}_{23}.
\end{gather}
\end{lema}
\begin{proof}
Let $e_{ab},e_{uv},e_{\alpha\beta}$ be elements of \B.\, We have to
prove
\begin{gather}
\tilde{\s}(e_{ab}\otimes m\langle
e_{uv},e_{\alpha\beta}\rangle)=(m\langle,\rangle\otimes
1)\tilde{\s}_{23}\tilde{\s}_{12},\label{ht1}\\
\tilde{\s}(m\langle e_{uv\otimes e_{\alpha\beta}}\otimes
e_{ab})=(1\otimes
m\langle,\rangle)\tilde{\s}_{12}\tilde{\s}_{23}(e_{uv}\otimes
e_{\alpha\beta}\otimes e_{ab}).\label{ht2}
\end{gather}
The left side of \eqref{ht1} is
\[
q^{c_{ab,\alpha v}+c_{ab,u\beta}}(q-q^{-1})e_{\alpha v}\otimes
e_{uv}\otimes e_{ab}
\]
while the right side is
\[
q^{c_{ab,uv}+c_{ab,\alpha\beta}}\langle e_{uv}\otimes
e_{\alpha\beta}\rangle\otimes e_{ab}
\]
Since \eqref{ces} we obtain $c_{ab,\alpha
v}+c_{ab,u\beta}=c_{ab,uv}+c_{ab,\alpha\beta}$. We conclude that
\eqref{ht1} holds. And by similar calculations \eqref{ht2} holds too.
\end{proof}
\begin{prop}\label{exten}
The linear morphism $\tilde{\s}$ can be extended to a linear morphism
\[
\s:\Uq\otimes_{k}\Uq\rightarrow\Uq\otimes_{k}\Uq
\]
where \Uq is the universal enveloping algebra of \slq, such that the
following diagram commutes,
\[\begin{CD}
\Uq\otimes_{k}\Uq @>\s>>\Uq\otimes_{k}\Uq\\
    @AAA                       @AAA\\
    \slq\otimes_{k}\slq @>\s>>\slq\otimes_{k}\slq
\end{CD}
\]
besides the multiplicativity conditions
\begin{gather}
\s(1\otimes m)=(m\otimes 1)(1\otimes\s)(\s\otimes 1)\label{m1}\\
\s(m\otimes 1)=(1\otimes m)(\s\otimes 1)(1\otimes\s)\label{m2}
\end{gather}
holds, where $m$ denotes the multiplication on \Uq.
\end{prop}

Now we may apply some remarks from Durdevi\'c \cite{Durdevic}. Suppose
that we are in the conditions of the proposition \ref{exten}. The
tensorial product $U(L)\otimes U(L)$ is a natural $U(L)$-bimodule and we
can define a multiplication over $U(L)\otimes U(L)$ by
\begin{equation}\label{mult}
(a\otimes b)(c\otimes d)=a\sigma(b\otimes c)d
\end{equation}
\begin{cor}
$\Uq\otimes_{k}\Uq$  is an   associative algebra with multiplication
defined by equation \eqref{mult}.
\end{cor}
\begin{proof}
The multiplicativity conditions \eqref{m1}, \eqref{m2} imply the
associativity of \eqref{mult}, \cite{Durdevic}
\end{proof}

In a similar way to the classical case define $\phi:\slq\rightarrow
\Uq\otimes_{k}\Uq$ by
\begin{equation}\label{coprod}
\phi(e_{i(i+1)})=e_{i(i+1)}\otimes 1+1\otimes e_{i(i+1)},\quad 1\leq
i\leq n.
\end{equation}
\begin{teo}
The morphism defined by \eqref{coprod}  can be extended to a morphism of
$k$-algebras
\[
\phi:\Uq\rightarrow \Uq\otimes_{k}\Uq
\]
\end{teo}
\begin{proof}
The $k$-algebra \Uq is generated by $x_{i}=e_{i(i+1)},$ $i=1,\ldots,n$
module the relations
\begin{gather}
x_{i}^{2}x_{j}-(q+q^{-1})x_{i}x_{j}x_{i}+x_{j}x_{i}^{2}=0\;\text{ if
}|i-j|=1,\\
x_{i}x_{j}-x_{j}x_{i}=0,\text{ if }|i-j|>1.
\end{gather}
We have to prove that $\phi$ preserves these relations. Let us put
$\sigma(x_{i}\otimes x_{j})=q^{c_{ij}}x_{j}\otimes x_{i},$ $1\leq
i,j,\leq n.$
\begin{multline*}
\phi(x_{i})^{2} \phi(x_{j})=\\
x_{i}^{2}x_{j}\otimes 1+x_{i}^{2}\otimes x_{j}+q^{c_{ij}}x_{i}x_{j}\otimes
x_{i}+x_{i}\otimes x_{i}x_{j}+q^{2+c_{ij}}x_{i}x_{j}\otimes
x_{i}+q^{2}x_{i}\otimes x_{i}x_{j}\\
+q^{2c_{ij}}x_{j}\otimes x_{i}^{2}+1\otimes x_{i}^{2}x_{j}
\end{multline*}
while
\begin{multline*}
\phi(x_{j})\phi(x_{i})^{2}=\\
x_{j}x_{i}^{2}\otimes 1+x_{j}x_{i}\otimes x_{i}+q^{2}x_{j}x_{i}\otimes
x_{i}+x_{j}\otimes x_{i}^{2}+q^{2cji}x_{i}^{2}\otimes
x_{j}+q^{2+c_{ji}}x_{i}\otimes x_{j}x_{i}\\
+q^{c_{ji}}x_{i}\otimes x_{j}x_{i}+1\otimes x_{j}x_{i}^{2}
\end{multline*}
and
\begin{multline}
\phi(x_{i})\phi(x_{j})\phi(x_{i})=\\
x_{i}x_{j}x_{i}\otimes 1+x_{i}x_{j}\otimes
x_{i}+q^{c_{ij}}x_{i}^{2}\otimes x_{j}+x_{i}\otimes x_{j}x_{i}\\
+q^{c_{ij}+2}x_{j}x_{i}\otimes x_{i}+q^{c_{ij}}x_{j}\otimes x_{i}^{2}\\
+q^{c_{ij}+2}x_{i}\otimes x_{i}x_{j}+1\otimes x_{i}x_{j}x_{i}
\end{multline}
It follows
\begin{multline}
\phi(x_{i})^{2}\phi(x_{j})+\phi(x_{j})\phi(x_{i})^{2}=\\
q^{-1}(q+q^{-1})(x_{i}^{2}\otimes x_{j})+(q^{-1}+q)x_{i}\otimes
x_{i}x_{j}+(q+q^{-1})x_{j}\otimes x_{i}^{2}\\
+(q+q^{-1})x_{i}x_{j}x_{i}\otimes 1 +1\otimes
(q+q^{-1})x_{i}x_{j}x_{i}\\
=(q+q^{-1})\phi(x_{i})\phi(x_{j})\phi(x_{i})
\end{multline}
\end{proof}

We shall define the counit. This can be done following the classical
case.

The commutative ring $k$ is a basic $T$-Lie algebra in the obvious way
and the zero morphism $0:\slq\rightarrow k$ is a morphism
$\epsilon=U(0):\Uq\rightarrow U(k)\simeq k$ of $k$-algebras.

\begin{define}
\[
x_{i}=e_{i(i+1)},\;i=1,\ldots,n.
\]
\end{define}

\begin{prop}
Let $\mathcal{C}=\{x_{1},\ldots,x_{n}\}.$ Then
\[
\phi(x_{i_{1}})\ldots\phi(x_{i_{m}})=x_{i_{m}}\ldots x_{i_{m}}\otimes 1
+\sum_{j}u_{j}\otimes v_{j}=1\otimes x_{i_{1}}\ldots
x_{i_{m}}+\sum_{l}a_{l}\otimes b_{l}
\]
where each $u_{j},v_{j},a_{l},b_{l}$ is a non-empty product of basic
elements in $\mathcal{C}$. It follows
\[
(1\otimes\epsilon)\phi(x_{i_{1}})\ldots\phi(x_{i_{m}})=x_{i_{1}}\ldots
x_{i_{m}}=(\epsilon\otimes 1)\phi(x_{i_{1}})\ldots \phi(x_{i_{m}}).
\]
Since $\mathcal{C}$ is a generator set of \Uq we get that $\epsilon$ is
the counit for $\phi$.
\end{prop}
\end{section}
\begin{section}{Antipode}
Let $L$ be \slq as a $T$-Lie algebra.
\begin{define}
\noindent
\begin{enumerate}
\item The {\it opposite} $T$-Lie algebra $L^{op}$ is defined as
\[
[,]^{op}=[,]_{q}S,\:S^{op}=S,\,\langle,\rangle^{op}=
\sigma^{-1}\langle,\rangle
\]
\item Let $m$ be the product of $U(L)$. The oppostive algebra
$U(L)^{op}$ is defined as $U(L)$ itself in its $k$-module structure
and product given by
\[
m^{op}=m\sigma.
\]
\end{enumerate}
\end{define}
\begin{prop}
Let $L$ be \slq.
\begin{enumerate}
\item The $k$-algebra $U(L)^{op}$ is associative and $L^{op}$ is a basic
$T$-Lie algebra.
\item The map $\eta:L\rightarrow L,\,x\mapsto -x$ is a $T$-Lie algebra
morphism.
\item There exist an isomorphism
\[
U(L^{op})\simeq U(L)^{op}
\]
of $k$-algebras.
\end{enumerate}
\end{prop}
\begin{prop}
There exist a morphism
\[
\eta:U(L)\rightarrow U(L)^{op}
\]
of $k$-algebras such that $\eta(x)=-x,\,\forall x\in L.$
\end{prop}

Now. let us consider the {\it quantum plane} $\mathbb{A}_{q}^{2|0}$ defined by
the ring
\[
\mathbb{A}_{q}^{2|0}=k\langle x,y/\langle yx-q xy\rangle
\]
where $k\langle x,y\rangle$ means an associative algebra freely
generated by $x,y$.

For positive integers $i\leq n$, we define the numbers
\[
\begin{pmatrix}
m\\
i
\end{pmatrix}_{q}
\]
by the equation in  $\mathbb{A}_{q}^{2|0}=k\langle x,y/\langle yx-q
xy\rangle,$
\[
(x+y)^{m}=\sum_{i=0}^{m}\begin{pmatrix}
                         m\\
                         i
                        \end{pmatrix}_{q}x^{m-i}y^{i}
\]
\begin{lema}
\begin{equation}\label{q-bin}\sum_{i=0}^{m}\begin{pmatrix}
                m\\
                i
              \end{pmatrix}_{q^{2}}(-1)^{i}q^{i(i-1)}=0
\end{equation}
\end{lema}
\begin{proof}
If $m=1$ the equation \eqref{q-bin} holds. Now, suppose \eqref{q-bin}.
Then,
\begin{multline}\label{inds}
\sum_{i=0}^{m+1}\begin{pmatrix}
                m+1\\
                i
                \end{pmatrix}_{q^{2}}(-1)^{i}q^{i(i-1)}=\\
\begin{pmatrix}
m+1\\
0
\end{pmatrix}_{q^{2}}+\sum_{i=1}^{m}\begin{pmatrix}
                                     m\\
                                     i
                                    \end{pmatrix}_{q^{2}}(-1)^{i}q^{i(i-1)}
+\begin{pmatrix}
m+1\\
m+1
\end{pmatrix}_{q^{2}}(-1)^{m+1}q^{(m+1)m}
\end{multline}
and because
\[
\begin{pmatrix}
m+1\\
i
\end{pmatrix}_{q}=\begin{pmatrix}
m\\i-1\end{pmatrix}_{q}+q^{i}\begin{pmatrix}m\\ i\end{pmatrix}_{q}\text{
for }1\leq i\leq m,
\]
(see \cite[p. 74]{Kassel})   then \eqref{inds} is equal to
\begin{align*}
1&+\sum_{i=1}^{m}\begin{pmatrix}m\\i-1\end{pmatrix}(-1)^{i}q^{i(i-1)}
+\sum_{i=1}^{m}\begin{pmatrix}m\\i\end{pmatrix}_{q^{2}}(-1)^{i}q^{i^{2}+i}+
(-1)^{n+1}q^{(n+1)n}\\
&=1+\sum_{j=0}^{m}\begin{pmatrix}m\\j\end{pmatrix}_{q^{2}}(-1)^{j+1}q^{(j+1)j}
+\sum_{i=1}^{m}\begin{pmatrix}m\\i\end{pmatrix}_{q^[2]}(-1)^{i}
q^{i(i+1)}+(-1)^{n+1}q^{n(n+1)}\\
&=0
\end{align*}
\end{proof}
\begin{lema}
If $\iota:U(L)^{op}\rightarrow U(L)$ is the natural $k$-module morphism
then
\begin{enumerate}
\item\[
\phi(x_{j})^{m}=\sum_{i=0}^{m}\begin{pmatrix}m\\i\end{pmatrix}_{q^{2}}
x_{j}^{m-i}\otimes x_{j}^{i},\;1\leq j\leq n
\]
\item\begin{multline*}
\xt{1}\xt{2}\ldots\xt{u}=\\
q^{\sum_{a<b}c_{j_{a}j_{b}}i_{a}(n_{j_{a}}-i_{b})}\xj{1}\xj{2}\ldots
\xj{u}\otimes \xij{1}\xij{2}\ldots\xij{u}
\end{multline*}
\item \begin{multline*}
(1\otimes\eta)(\xj{1}\xj{2}\ldots\xj{u}\otimes
\xij{1}\xij{2}\ldots\xij{u})=\\
q^{\sum_{a<b}c_{j_{a}j_{b}}i_{a}i_{b}}\xj{1}\xj{2}\ldots\xj{u}\otimes\iota
\eta(\xij{u})\ldots\iota\eta(\xij{1})
\end{multline*}
\item For $j=1,\ldots,n,$
\[
\eta(x_{j}^{i})=(-1)^{i}q^{i(i-1)}(\iota x_{j})^{i}.
\]
\item Denote with $m$ the product of \Uq. Then
\[
m(1\otimes\eta)\phi(\xij{1}\ldots\xij{u})=0,\;1\leq j_{1},\ldots,
j_{u}\leq n.
\]
if $n_{1},\ldots,n_{u}$ are all positive integers.
\end{enumerate}
\end{lema}
\begin{proof}
By straightforward computations. By example, if
$c=\sum_{a<b}c_{j_{a}j_{b}i_{a}(n_{j_{a}}-i_{b})}$ then
\begin{multline}
m(1\otimes\eta)\phi(\xij{1}\ldots\xij{u})=m(1\otimes\eta)\phi(x_{j_{1}})^{n_1}
\ldots\phi(x_{j_{u}})^{n_{u}}=\\
\sum_{i_{1},\ldots i_{u}=0}^{n_{1},\ldots,n_{u}}\qbin{1}\ldots\qbin{u}
q^{c}m(\xj{1}\ldots)\xj{u}\otimes \eta(\xij{u}\ldots\xij{1})\\
=\sum_{i_{1},\ldots i_{u}=0}^{n_{1},\ldots,n_{u}}\qbin{1}\ldots\qbin{u}
q^{\sum_{a<b}c_{j_{a}j_{b}}i_{a}n_{j_{b}}}\xj{1}\ldots\xj{u}\iota\eta(\xij{u})\ldots
\iota\eta(\xij{1}) \\
=\sum_{i_{1},\ldots i_{u-1}=0}^{n_{1},\ldots,n_{u-1}}\qbin{1}\ldots\qbin{u}
q^{\sum_{a<b}c_{j_{a}j_{b}}i_{a}n_{j_{b}}}
\xj{1}\ldots\xj{u-1}\\
\sum_{i_{u}=0}^{n_{u}}\qbin{u}(-1)^{i_{u}}q^{i_{u}(i_{u}-1)}\iota x_{j_{u}}^{n_{u}}
\iota\eta(\xij{u-1})\ldots
\iota\eta(\xij{1})=0
\end{multline}
\end{proof}
\begin{prop}
Let $\iota:U(L)^{op}\rightarrow U(L)$ the natural $k$-morphism. Then,
the $k$-module morphism
\[
\kappa=\iota\eta:U(L)\rightarrow U(L)
\]
is the antipode for the coproduct $\phi$.
\end{prop}
\begin{proof}
\end{proof}
\end{section}
\begin{section}{The additional condition of the braided quantum group
definition}
\begin{lema}\label{pread}
If
\begin{align}
(\sigma\otimes 1)(1\otimes \sigma)(\phi\otimes
1)=(1\otimes\phi)\sigma;\\
(1\otimes\sigma)(\sigma\otimes 1)(1\otimes\phi)=(\phi\otimes 1)\sigma.
\end{align}
then
\begin{equation}\label{ad}
(\sigma\otimes 1^{2})(1\otimes\phi\otimes 1)(\sigma^{-1}\otimes
1)(1\otimes\phi)=(1^{2}\otimes\sigma)(1\otimes\phi \otimes
1)(1\otimes\sigma^{-1})(\phi\otimes 1)
\end{equation}
\end{lema}
\begin{proof}
\begin{align*}
(\sigma\otimes 1^{2})(1\otimes\phi\otimes 1)(\sigma^{-1}\otimes
1)(1\otimes\phi)& =(\,(\sigma\otimes 1)(1\otimes\phi)\otimes
1\,)(\sigma^{-1}\otimes 1)(1\otimes\phi)\\
&= (\,(1\otimes\sigma^{-1})(\phi\otimes 1)\sigma\otimes
1\,)(\sigma^{-1}\otimes 1)(1\otimes\phi)\\
&=(1\otimes\sigma^{-1}\otimes 1)(\phi\otimes\phi),
\end{align*}
on the other hand,
\begin{align*}
(1^{2}\otimes\sigma)(1\otimes\phi \otimes
1)(1\otimes\sigma^{-1})(\phi\otimes
1)&=(\,1\otimes(1\otimes\sigma)(\phi\otimes
1)\,)(1\otimes\sigma^{-1})(\phi\otimes 1)\\
&=(\,1\otimes(\sigma^{-1}\otimes
1)(1\otimes\phi)\sigma\,)(1\otimes\sigma^{-1})(\phi\otimes 1)\\
&=(1\otimes\sigma^{-1}\otimes 1)(\phi\otimes\phi).
\end{align*}
We conclude that \eqref{ad} holds.
\end{proof}
\begin{prop}
The following equations holds,
\begin{gather}
(\sigma\otimes 1)(1\otimes\sigma)(\phi\otimes 1)=(1\otimes\phi)\sigma\\
(1\otimes\sigma)(\sigma\otimes 1)(1\otimes\phi)=(\phi\otimes
1)\sigma.\label{too}
\end{gather}
\end{prop}
\begin{proof}
Suppose that $x_{j_[1]},\ldots,x_{j_{u}}$  and
$x_{k_{1}},\ldots,x_{k_{v}}$ are elements of $\{x_{1},\ldots,x_{n}\}.$
Put $x_{J}=\xnj{1}\ldots\xnj{u}$, $x_{K}=\xmk{1}\ldots\xmk{v}$. Then
\begin{multline}
(\phi\otimes 1)(x_{J}\otimes x_{K})=\sum_{i_{1}=0,\ldots
i_{u}=0}\qbin{1}\ldots\qbin{u}q^{\sum_{a<b}c_{j_{a}j_{b}}i_{a}(n_{j_{b}}
-i_{b})}\\
\xj{1}\ldots\xj{u}\otimes \xij{1}\ldots\xij{u}\otimes
\xmk{1}\ldots\xmk{v}
\end{multline}
it follows,
\begin{multline*}
(1\otimes\sigma)(\phi\otimes 1)(x_{J}\otimes x_{K})=
\sum_{i_{1}=0,\ldots
i_{u}=0}\qbin{1}\ldots\qbin{u}q^{\sum_{a<b}c_{j_{a}j_{b}}i_{a}(n_{j_{b}}
-i_{b})}\\
q^{\sum_{a,b}c_{j_{a}k_{b}}i_{a}m_{b}}\xj{1}\ldots\xj{u}\otimes
\xmk{1}\ldots\xmk{v}\otimes \xij{1}\ldots\xij{u}
\end{multline*}
and
\begin{multline*}
(\sigma\otimes 1)(1\otimes\sigma)(\phi\otimes 1)(x_{J}\otimes x_{K})=\\
\sum_{i_{1}=0,\ldots
i_{u}=0}\qbin{1}\ldots\qbin{u}q^{\sum_{a<b}c_{j_{a}j_{b}}i_{a}(n_{j_{b}}
-i_{b})}q^{\sum_{a,b}(c_{j_{a}k_{b}}i_{a}m_{b}
+c_{j_{a}k_{b}}(n_{a}-i_{a})m_{b})} \\
\xmk{1}\ldots\xmk{v}\otimes\xij{1}\ldots\xij{u}\otimes\xj{1}\ldots\xj{u}\\
=q^{\sum_{a,b}c_{j_{a}k_{b}}n_{a}m_{b}}\xmk{1}\ldots\xmk{v}\otimes
\phi(x_{j_{1}}^{n_{1}})\ldots\phi(x_{j_{u}}^{n_{u}})\\
=(1\otimes\phi)\sigma(x_{J}\otimes x_{K})
\end{multline*}
By similar calculations \eqref{too} hold too.
\end{proof}

Suppose that $A,B$ are $k$-algebras (non-associative, perhaps), and
$\beta:A\otimes_{k}B\rightarrow B\otimes_{k}A$ a $k$-morphism. Moreover,
denote with $A\otimes_{\beta}B$ to $A\otimes_{k}B$ in its $k$-module
structure and product given by
\[
(a\otimes b)(c\otimes d)=a\beta(b\otimes c)d
\]
\begin{teo}
The algebra \Uq is a braided quantum group with
\begin{enumerate}
\item coproduct:
\[
\phi:\Uq\rightarrow\Uq\otimes_{\sigma}\Uq
\]
induced by $\phi(x_{i})=x_{i}\otimes 1+1\otimes x_{i},\;i=1,\ldots,n;$
\item counit:
\[
\epsilon:\Uq\rightarrow k
\]
induced by $\epsilon(x)=0,\;\forall x\in\slq$;
\item antipode:
\[
\eta:\Uq\rightarrow\Uq
\]
induced by $\eta(x)=-x,\;\forall x\in\slq.$
\end{enumerate}
\end{teo}
\end{section}
\bigskip

\begin{center}
{\Large\bf Acknowlegment}
\medskip
\end{center}

I would like to thank M. Durdevi\'c. The existence of a braided
product \eqref{mult} in the case $\mathbb{A}_{n}$ was suggested by Durdevi\'c,
besides the lemma \ref{pread} belongs to him.


\begin{thebibliography}{9}
\bibitem{yo} C. Bautista, A Poincar\'e-Birkhoff-Witt theorem for
generalized color Lie algebras, {\it preprint} q-alg/9706016. Submitted to the
{\it J. of Math. Phys.}
\bibitem{Bourbaki}N. Bourbaki, ``Lie Groups and Lie Algebras,'' Elements
of Mathematics. Chap. 1-3. Springer-Verlag. Great Britain. 1989.
\bibitem{Durdevic}M. Durdevi\'{c}, On braided Quantum Groups, {\it
preprint} q-alg/9412003.
\bibitem{Gurevich} D. Gurevich, The Yang-Baxter Equation and a
Generalization of Formal Lie Theory, {\it Soviet Math. Dokl.} {\bf 33} (1986),
No. 758-762.
\bibitem{Kassel} C. Kassel, ``Quantum Groups,'' GTM 155. Springer-Verlag.
New York.  1995.
\bibitem{Street}R. Street. ``Quantum Groups: an entre\'e to modern
algebra,'' Preliminary version. Macquarie University. 1993-1994.
\end{thebibliography}
\end{document}